\documentclass[11 pt]{amsart}
\usepackage[applemac]{inputenc}
\RequirePackage{etex}

\usepackage[dvipsnames,svgnames,x11names,hyperref]{xcolor}
\usepackage[colorlinks=true,linkcolor=Maroon,citecolor=blue,urlcolor=blue,hypertexnames=false,linktocpage]{hyperref}
\usepackage{bookmark}
\usepackage{amsthm,thmtools,amssymb,amsmath,amscd}

\usepackage{fancyhdr}
\usepackage{esint}
\usepackage{enumerate}

\usepackage{pictexwd,dcpic}
\usepackage{graphicx}
\usepackage{caption}
\swapnumbers
\declaretheorem[name=Theorem,numberwithin=section]{thm}
\declaretheorem[name=Remark,style=remark,sibling=thm]{rem}
\declaretheorem[name=Lemma,sibling=thm]{lemma}
\declaretheorem[name=Proposition,sibling=thm]{prop}
\declaretheorem[name=Definition,style=definition,sibling=thm]{defn}
\declaretheorem[name=Corollary,sibling=thm]{cor}

\numberwithin{equation}{section}

\usepackage{cleveref}
\crefname{lemma}{Lemma}{Lemmata}
\crefname{prop}{Proposition}{Propositions}
\crefname{thm}{Theorem}{Theorems}
\crefname{cor}{Corollary}{Corollaries}
\crefname{defn}{Definition}{Definitions}
\crefname{example}{Example}{Examples}
\crefname{rem}{Remark}{Remarks}
\crefname{assum}{Assumption}{Assumptions}
\crefname{nota}{Notation}{Notation}

\usepackage{autonum}

\newcommand{\ti}{\tilde}

\newcommand{\cn}{\colon}
\newcommand{\sub}{\subset}

\newcommand{\bbR}{\mathbb{R}}

\newcommand{\bbS}{\mathbb{S}}

\newcommand{\8}{\infty}

\newcommand{\al}{\alpha}
\newcommand{\be}{\beta}
\newcommand{\ga}{\gamma}
\newcommand{\de}{\delta}
\newcommand{\ep}{\epsilon}
\newcommand{\ka}{\kappa}
\newcommand{\la}{\lambda}
\newcommand{\om}{\omega}
\newcommand{\si}{\sigma}
\newcommand{\Si}{\Sigma}

\newcommand{\Om}{\Omega}
\newcommand{\De}{\Delta}
\newcommand{\Ga}{\Gamma}

\newcommand{\La}{\Lambda}

\newcommand{\cN}{\mathcal{N}}
\newcommand{\cU}{\mathcal{U}}

\newcommand{\cO}{\mathcal{O}}

\newcommand{\cP}{\mathcal{P}}
\newcommand{\cC}{\mathcal{C}}
\newcommand{\cX}{\mathcal{X}}

\newcommand{\del}{\partial}
\newcommand{\n}{\nabla}
\newcommand{\II}[2]{\mrm{II}\br{#1,#2}}
\newcommand{\fa}{\forall}
\newcommand{\rt}{\sqrt}

\newcommand{\ip}[2]{\langle #1,#2 \rangle}
\newcommand{\fr}[2]{\frac{#1}{#2}}
\newcommand{\tfr}[2]{\tfrac{#1}{#2}}
\newcommand{\x}{\times}

\DeclareMathOperator{\dive}{div}

\DeclareMathOperator{\pr}{pr}

\DeclareMathOperator{\tr}{tr}

\DeclareMathOperator{\Rc}{Rc}

\newcommand{\pf}[1]{\begin{proof}#1 \end{proof}}
\newcommand{\eq}[1]{\begin{equation}\begin{alignedat}{2} #1 \end{alignedat}\end{equation}}

\newcommand{\br}[1]{\left(#1\right)}
\newcommand{\abs}[1]{\lvert #1\rvert}
\newcommand{\enum}[1]{\begin{enumerate}[(i)] #1 \end{enumerate}}

\newcommand{\ra}{\rightarrow}
\newcommand{\hra}{\hookrightarrow}

\newcommand{\mt}{\mapsto}

\newcommand{\mc}{\mathcal}

\newcommand{\mrm}{\mathrm}

\newcommand{\hp}{\hphantom}
\newcommand{\q}{\quad}

\usepackage{slashed}
\usepackage{accents}
\newcommand{\ubar}[1]{\underaccent{\bar}{#1}}
\newcommand{\s}{\slashed}
\newcommand{\sn}{\s\n}

\begin{document}
\title[Mean curvature flow in null hypersurfaces]{Mean curvature flow in null hypersurfaces and the detection of MOTS}
\author{Henri Roesch}
\address{Department of Mathematics, Columbia University, Broadway 2990, New York, NY 10027, USA}
\email{{\href{mailto:roesch@math.columbia.edu}{roesch@math.columbia.edu}}}
\author{Julian Scheuer}
\address{School of Mathematics, Cardiff University, Senghennydd Road, Cardiff CF24 4AG, Wales}
\email{{\href{mailto:scheuerj@cardiff.ac.uk}{scheuerj@cardiff.ac.uk}}}
\date{\today}
\thanks{Funded by the ``Deutsche Forschungsgemeinschaft" (DFG, German research foundation); Project "Quermassintegral preserving local curvature flows"; No. SCHE 1879/3-1.}
\thanks{Funded by the National Science Foundation under award DMS-1703184.}

\begin{abstract}
We study the mean curvature flow in 3-dimensional null hypersurfaces. In a spacetime a hypersurface is called null, if its induced metric is degenerate. The speed of the mean curvature flow of spacelike surfaces in a null hypersurface is the projection of the codimension-two mean curvature vector onto the null hypersurface. We impose fairly mild conditions on the null hypersurface. Then for an outer un-trapped initial surface, a condition which resembles the mean-convexity of a surface in Euclidean space, we prove that the mean curvature flow exists for all times and converges smoothly to a marginally outer trapped surface (MOTS). As an application we obtain the existence of a global foliation of the past of an outermost MOTS, provided the null hypersurface admits an un-trapped foliation asymptotically.      
\end{abstract}

\maketitle

\section{Introduction}

Let $(M,g)$ be a four dimensional, time-oriented Lorentzian manifold or \textit{spacetime} with Levi-Civita connection $D$, where for convenience we write 
\eq{\langle X,Y\rangle:=g(X,Y)} for vector fields $X,Y$ of $M$. Let $\Sigma\subset M $ be an embedded spacelike 2-sphere. The second fundamental form of $\Sigma$ in $M$ and the corresponding mean curvature vector are denoted by
\eq{\II{V}{W} = (D_{V}W)^{\perp},\quad \vec{H} = \tr_{\Sigma}\mathrm{II},}
for all sections $V,W\in \Gamma(T\Sigma)$. 
For a future directed null normal $l\in\Gamma(T^\perp\Sigma)$, 
the conditions 
\eq{-\langle l,\vec{H}\rangle > 0, \q-\langle l,\vec{H}\rangle<0}on $\Si$ are called \textit{outer un-trapped}, or \textit{outer trapped} respectively. Those properties are referred to hold {\it{weakly}}, provided the respective weak inequality holds. If the space-time $(M,g)$ is globally hyperbolic satisfying the null energy condition, then the famous singularity theorems of Hawking and Penrose imply $(M,g)$ forms a singularity in the causal future of an outer trapped surface $\Sigma$. A \textit{marginally outer trapped surface}, or MOTS, is identified by the condition 
\eq{\langle l,\vec{H}\rangle = 0.} 
If a MOTS bounds a family of trapped surfaces and $(M,g)$ is asymptotically predictable, one concludes that a MOTS lies within the black-hole region of $(M,g)$. Consequently, MOTSs are studied as `quasi-local versions' of the event horizon of $(M,g)$, the boundary of the black-hole region.

Our main motivation is to use a geometric flow in locating a MOTS. This idea goes back to work of Tod \cite{Tod:/1991}, where he suggested the mean curvature flow to find a MOTS inside a time-symmetric slice of spacetime. Time-symmetry indicates a Riemannian hypersurface fully decoupled (or totally geodesic) from the ambient spacetime geometry. Within such a slice, a MOTS is identified as a minimal surface. Given an initial boundary on which to initiate mean curvature flow, work by White \cite{White:/2000} showed, provided this region encloses a minimal surface, that an outermost minimal surface will result in the limit of the flow. In the non time-symmetric case, whereby the Riemannian slice observes a non-trivial second fundamental form in spacetime, Tod also suggested the use of \textit{null mean curvature flow}. Bourni-Moore \cite{BourniMoore:/2019} subsequently formulated a theory of weak solutions to the null mean curvature flow in this setting. Their approach similarly used a weak level set formulation as in the famous work of Huisken-Ilmanen \cite{HuiskenIlmanen:/2001} (in the time-symmetric case), and then Moore \cite{Moore:/2014} (in the non-time symmetric case) to address inverse mean curvature flow. Governed by a scalar degenerate elliptic equation, Bourni-Moore were able to use elliptic regularization to establish existence of a weak solution to the null mean curvature flow. They observe convergence to a measure theoretic `generalized MOTS' lying outside the outermost MOTS. Their solution also exhibits a blow-up at the outermost MOTS due to its strong association with a solution to Jang's equation. In turn, the blow-up of solutions to Jang's equation is a key property in characterizing a MOTS as in the celebrated proof of the Positive Mass Theorem by Schoen-Yau \cite{SchoenYau:/1979, SchoenYau:/1981}. We also mention the work of Eichmair \cite{Eichmair:/2009} where this blow-up was used to solve the Plateau problem for MOTSs. More specifically, one considers a region within a non time-symmetric slice bounded by an outer trapped region at one end and an outer un-trapped region at the other. Eichmair then developed a versatile technique to force and control a blow-up of Jang's equation, subsequently showing the existence of a MOTS within this region. 

In this paper we propose a new method to find MOTSs, namely by employing the mean curvature flow (MCF) in a null hypersurface of a spacetime. A null hypersurface is characterized by the property that the induced metric inherited from the spacetime is degenerate and
more specifically, we will consider null hypersurfaces, which are foliated by spherical leaves. 
MCF has been widely studied as a flow of hypersurfaces in Riemannian and Lorentzian manifolds. If $x$ denotes a time-dependent family of embeddings of a smooth manifold into a Riemannian or Lorentzian ambient space, MCF can concisely be written as
\eq{\dot x = \De x = \vec{H},}
where $\De$ is the Laplace-Beltrami operator with respect to the metric induced by $x(t,\cdot)$ and $\dot{x} = \dot{x}(\cdot,\xi)$ is the velocity of the curve $t\mt x(t,\xi)$, $\xi$ being an element of the embeddings' common domain, which will be $\bbS^{2}$ in this paper. For any embedding, $\De x$ is perpendicular to the hypersurface and hence proportional to a chosen normal vector field. Then
\eq{\label{MCF-nondeg}\dot x = -\si H\nu,}
where $\si = \ip{\nu}{\nu}$ is the signature of the ambient space and $H = \abs{\vec{H}}.$
This is one of the most important equations of geometric analysis and the literature is vast and exponentially growing. We do not give a very detailed account here. In the smooth setting it was pioneered by Huisken for convex hypersurfaces of the Euclidean space \cite{Huisken:/1984} and for entire graphs by Ecker/Huisken \cite{EckerHuisken:11/1989}. Huisken has also developed surgery in the $2$-convex setting with Sinestrari \cite{HuiskenSinestrari:/2009}. Other important aspects of MCF relate to geometric pinching estimates \cite{Andrews:/2012,Brendle:12/2016,Langford:08/2017}, Harnack inequalities \cite{Hamilton:/1995} and ancient solutions \cite{BourniLangfordTinaglia:07/2020}. 

 In the Lorentzian setting there are convergence results for spacelike entire graphs \cite{Ecker:/1997,LambertLotay:02/2021}. In higher codimension the picture is much less developed and good results are usually only available under strong pinching conditions on the initial hypersurface, e.g. \cite{AndrewsBaker:/2010,BakerNguyen:11/2017}.
 
To the best of our knowledge, MCF has never been studied as a flow within a null hypersurface $\cN$ as we propose to do it here. The problems are obvious: The normal to any spacelike surface within $\cN$ is a null vector and hence a representation of MCF in either of the above forms is impossible. Neither is there a Levi-Civita connection nor a unit normal vector field. We believe that the best way to write MCF in this setting is to take
\eq{\label{MCFN}\dot x = \pr_{T\cN}\vec H,}
where $\pr_{T\cN}$ is the skew-orthogonal projection of a vector onto $T\cN$, see \Cref{sec:setup}. About this approach there are good and bad news. The bad news are of geometric nature: We pick up classical ``higher-codimension-problems'', such as the presence of torsion, that we have to deal with. The good news are of PDE-nature: Spacelike MCF in our null hypersurface is automatically graphical and the equation does not see the {\it{slope}} of the graphs, because the flow direction is a null vector. This makes things easier from a PDE point of view. To explain this further, note that graphical MCF in Euclidean (Lorentzian) space can be written as
\eq{\del_{t}\om = -\si H\rt{1\pm\abs{\n\om}^{2}},}
where $\om$ is the graph function of the flow hypersurfaces. However, as we shall see later in \Cref{sec:setup}, in a suitable gauge MCF in a null hypersurface is given by
\eq{\del_{t}\om = -H.}
It is interesting to see how this flow somehow seems to ``interpolate'' between its Riemannian and Lorentzian relatives.
It is important to note that this flow differs entirely from the previously mentioned null mean curvature flow by Bourni/Moore \cite{BourniMoore:/2019}, as their flow is a variation of hypersurfaces in a Riemannian manifold.  

 In this paper we show that \eqref{MCFN} is capable of doing the following:
Given a null hypersurface $\cN$ supporting a trapped surface or MOTS, we identify fairly generic constraints on $\cN$ for which our mean curvature flow \eqref{MCFN} from any outer un-trapped initial cross-section exists for all times and converges smoothly to a MOTS. 

\subsection*{Main results}
We prove the long-time existence of the mean curvature flow for spacelike spherical cross-sections within a null hypersurface $\cN$ and show that it converges to a MOTS. To formulate the result, we go over notation very briefly. For a detailed description see \Cref{sec:setup}. In the following theorem, $R$ denotes the Riemann tensor as defined in \eqref{Rm} and $G$ is the associated Einstein tensor. The past directed null vector $k$ is part of a {\it{null basis}} $\{k,l\}\sub\Gamma(T^{\perp}\Sigma)$, $\Si\cong\bbS^2$, given by:
\eq{\label{null partner}\ip{l}{l} = 0,\q \ip{k}{l}=2,\q \ip{k}{k}=0.} If we also denote by $k$ the unique null geodesic vector field extension throughout $\cN$, then we may rescale to a vector field $\ubar L:= ak$, for some $a\in C^\infty(\cN)$, $a>0$. The 2-tensor $\ubar\chi$ represents the second fundamental form of the null hypersurface with respect to $\ubar L$. We take $\hat{\ubar\chi}$ to be its traceless part. Finally,
\eq{\ubar\al(V,W)=\ip{R_{\ubar L V}\ubar L}{W}.}
Here is our main result.

\begin{thm}\label{thm:main}
Let $(M,g)$ be a 4-dimensional, time-oriented Lorentzian manifold, $\Si_{0}\sub M$ a weakly outer trapped two-sphere with respect to a future directed null normal section $l$, and let $\cN$ be the null hypersurface generated by the past directed null partner $k$ of $l$. Now consider $\ubar L = ak$, for $a\in C^\infty(\cN)$, $a>0$, satisfying the gauge condition:
\eq{\label{gauge}G(\ubar L,\ubar L)-d(2\ka-\tr\ubar\chi)(\ubar L)\geq \abs{(\tr\ubar\chi-4\ka)\hat{\ubar\chi}}+2\abs{\hat{\ubar\al}}+\tfr52\abs{\hat{\ubar\chi}}^{2},}
whereby $\ka = da(k)$.
 Then, if the null hypersurface $\Om\sub\cN$ generated by $\ubar L$ and $\Si_{0}$ admits an outer un-trapped cross-section $\Si_{\om_{0}},$ the mean curvature flow
\eq{\dot{x}\label{NHMCF} = \tfr 12 \ip{\vec H}{L_{\om_{t}}}\ubar L}
 initiated at $\Si_{\om_{0}}\sub\Om$ exists for all times and converges smoothly to a MOTS.
\end{thm}

The subscript $\om_{0}$ to the initial hypersurface indicates that $\Si_{\om_{0}}$ is given as the graph of a function $\om_{0}$ on $\Si_{0}$, a property that is preserved throughout the flow. The flow hypersurfaces are then denoted by $\Si_{\om_{t}}$, while $L_{\om_{t}}$ is the unique null partner of $\ubar L$ with respect to $\Si_{\om_{t}}$ as in \eqref{null partner}.
 All tensor norms $\abs{\cdot}$ in this and the following theorem are with respect to the induced metric on the background foliation of $\cN$, see \Cref{sec:setup} for a detailed account.

The gauge condition \eqref{gauge} described in \Cref{thm:main} translates to a family of ODE inequalities on the scaling $a\in C^\infty(\cN)$ along $k$ for which a wealth of solutions exist. Associated to each solution is a neighborhood of $\Si_0$ in $\cN$ and the $\ubar L$-directed one-sided part of that neighborhood we call $\Om$. In other-words, \Cref{thm:main} indicates that the mean curvature flow exists and converges to a MOTS provided an outer un-trapped cross-section exists `within reach of $\Si_0$' by one of these neighborhoods.\\
\indent For the 2-tensor $K$, representing the second fundamental form of $\cN$ with respect to $k$, we can also observe a sufficient condition satisfying the gauge constraint of \Cref{thm:main} as an energy type condition on specialized null structures. More specifically, in the case that $\cN$ is a Null Cone, whereby $\tr K>0$ throughout $\cN$: 

\begin{thm}\label{thm:cone}
	Let $\cN$ be a Null Cone admitting a weakly outer trapped cross-section $\Si_0\sub\cN$, and an outer un-trapped cross-section $\Si_{\om_{0}}\sub\cN$ to the timelike past of $\Si_0$. If the region bounded by $\Si_0,\Si_{\om_{0}}$ satisfies
	\eq{\label{energy}G(k,k)\geq \tfr52\abs{\hat{K}}^2+2\tr K\abs{\hat{K}}+2\abs{\hat\al},}
whereby $\al(V,W) := \ip{R_{kV}k}{W}$, then the mean curvature flow \eqref{NHMCF} initiated at $\Si_{\om_{0}}\sub\cN$ exists for all times and converges smoothly to a MOTS.
\end{thm}

\begin{rem}
\enum{
\item Note, in Schwarzschild space we have $G = \hat K = \hat\al = 0$.
\item
The assumption of $\Si_{0}$ being a topological two-sphere is not essential. We include this assumption because we formally rely on several calculations from \cite{Roesch:/2017}, where it is a standing assumption throughout the paper. 
For convenience of the reader, we do not want to make any statements which can not easily be checked from the given references, so we stick to this assumption. However, we believe it is not necessary.
}
\end{rem}

We also obtain a past foliation of an outermost MOTS, provided $\cN$ admits an asymptotically un-trapped foliation to the past of the MOTS. We refer to \Cref{sec:foliation} for a details. Briefly, we assume $\cN\cong (\La_-,\La_+)\times\bbS^2$, where the interval $(\La_-,\La_+)$ is generated by the vector field $k$ ($\La_+=\infty$ included).

\begin{thm}\label{thm:foliation}
Let $M$ and $\cN$ satisfy the conditions of \Cref{thm:main}, further we assume that $\cN$ admits a outer un-trapped foliation in a neighborhood of ${\La}_+$. Then there exists a strict outermost MOTS $\Si_{out}\sub\cN$ and a global foliation of $\cN$ by outer un-trapped surfaces to the timelike past of $\Si_{out}$.
\end{thm}

The paper is organized as follows. In \Cref{sec:setup} we introduce the underlying null geometry and the mean curvature flow equation in a detailed manner. In \Cref{Estimates} we prove the required a priori estimates, which are estimates up to $C^{1}$ in view of the quasi-linearity of the equation. In \Cref{proof main} we complete the proof of \Cref{thm:main}, in \Cref{sec:cone} we prove \Cref{thm:cone}, while in \Cref{sec:foliation} we prove \Cref{thm:foliation}.

\section{Setup}\label{sec:setup}

Let $(M,g)$ be a four dimensional, time-oriented Lorentzian manifold with Levi-Civita connection $D$, the convention from \cite{ONeill:/1983} for the Riemann tensor
\eq{\label{Rm}R_{XY}Z:= D_{[X,Y]}Z-[D_X,D_Y]Z}
and Einstein tensor 
\eq{G = \Rc - \tfr 12 S g,}
where 
\eq{\Rc(X,Y) = \tr_{g}\ip{R_{X\cdot}Y}{\cdot},\q S = \tr_{g}\Rc.} 
 Let $\Si\hra M$
be the embedding of a spacelike two-sphere. The second fundamental form of $\Si$ in $M$ and the corresponding mean curvature vector are denoted by
\eq{\II{V}{W}=(D_{V}W)^{\perp},\q \vec{H}=\tr_{\Si}\mrm{II},}
for all $V,W\in \Ga(T\Si)$. We mostly follow the notation from \cite{Roesch:/2017}.

\subsection*{Null Geometry} We now briefly describe null hypersurfaces in $(M,g)$. By definition, a null hypersurface $\cN \sub M$ is a smooth hypersurface such that the induced metric $g|_{\cN}$ is degenerate. We also assume $\cN$ is orientable. We may therefore observe a global vector field $\ubar L\in \Ga(T\cN)$, such that 
\eq{\ubar L_p^\perp = T_p\cN} for each $p\in \cN$. In particular this gives $\langle\ubar L,\ubar L\rangle = 0$. For $p\in \cN$, the hypersurface structure of $\cN$ also ensures a neighborhood $\cU$ and a smooth function $\upsilon\in C^\infty(\cU)$ such that 
\eq{\mathcal{V}:=\cU\cap\cN = \{\upsilon=0\}.} Moreover, if $q\in\mathcal{V}$, we have $(D\upsilon)_q^\perp = T_q\cN$, where we denote $D\upsilon:=\text{grad}(\upsilon)$. Consequently, the non-degeneracy of the ambient metric $g$ enforces that 
\eq{(D\upsilon)_q\propto \ubar L_q,}
 giving $\langle D\upsilon,D\upsilon\rangle|_\mathcal{V} \equiv 0$. From the identity
\eq{D_{D\upsilon}D\upsilon = \tfrac12D\langle D\upsilon,D\upsilon\rangle}
we conclude therefore that $T_q\cN\sub(D_{D\upsilon}D\upsilon)^\perp_q$ giving 
\eq{(D_{D\upsilon}D\upsilon)_q\propto (D\upsilon)_q,\q\text{equivalently}\q D_{\ubar L}\ubar L = \kappa \ubar L,}
for some $\kappa\in C^\infty(\cN)$. It follows that integral curves along $\ubar L$ are pre-geodesic, and $\cN$ is in-fact ruled by null geodesics. In the theory of general relativity, $\cN$ represents the geometry of a given light-ray congruence in spacetime.\\
\indent Now we introduce the \textit{second fundamental form} of $\cN$, for $X,Y\in\Ga(T\cN)$:
\eq{\ubar\chi(X,Y):=\langle D_X\ubar L,Y\rangle.}
This symmetric 2-tensor is defined up-to a scaling of $\ubar L$. For $X,Y\in \Ga(T\cN)$, $c\in C^\infty(\cN)$, we also observe the properties:
    \eq{\langle X+c\ubar L,Y\rangle = \langle X,Y\rangle,\,\,\,
    \ubar\chi(X+c\ubar L,Y) = \ubar\chi(X,Y).}
It follows that both the induced metric and the second fundamental form at $p\in\cN$ are fully characterized `modulo $\ubar L$'. Equivalently, both the metric and second fundamental form on $\cN$ at a point, $p\in\cN$, are fully determined by their restrictions to any spacelike slice $\Si$ through $p$. Consequently, whenever convenient we use the rather small abuse of notation to denote both the induced metric of any spacelike submanifold $\Si\sub \cN$ at $p$ and $(g|_\cN)_p$, by $\ga_p$. Similarly, we will simply denote by $\ubar\chi_p$ the restriction of $\ubar\chi$ to any spacelike $\Si\sub\cN$ at $p$. We also notice the function 
\eq{p\mt\tr_\Si\ubar\chi(p)} is independent of any spacelike slice $\Si$ through $p$, giving a well defined function $\tr\ubar\chi \in C^\infty(\cN)$. Finally, we bring to the attention of the reader that due to the symmetries of the Riemann curvature tensor, the same notational conventions (as for $\ga,\ubar\chi$) may be adopted for the 2-tensor 
\eq{\ubar\al(V,W) = \ip{R_{\ubar L V}\ubar L}{W}.}

\subsection*{The background foliation}
We will now assume a similar construction as in \cite{MarsSoria:09/2015}. We may choose a past-pointing null geodesic vector field $k\in\Ga(T\cN)$, such that $D_kk = 0$. We then assume the existence of a spacelike 2-sphere $\Si_0\sub\cN$ and the property that any geodesic along $k$ intersects $\Si_0$ precisely once. Consequently, $\cN$ is ruled by geodesics of $k$, denoted $\be_p(\la)$ for $p\in\Si_0$, whereby $\be_p(0) = p$. We also observe 
\eq{{\La}_+:=\inf\{{\la}_+(p)|p\in\Si_0\}>0,} whereby ${\la}_+(p):=\sup\{\la|\be_p(\la)\in\cN\}$. Similarly, 
\eq{{\La}_-:=\sup\{{\la}_-(p)|p\in\Si_0\}<0,} whereby ${\la}_-(p):=\inf\{\la|\be_p(\la)\in \cN\}$. Standard ODE theory ensures the mapping 
\eq{\label{null diffeo}({\La}_-,{\La}_+)\times\Si_0\to\cN,} given by $(\la,p)\to\beta_p(\la)$, is a smooth embedding onto an open subset of $\cN$. This open subset depends on both our choice of $\Si_0$ and our geodesic generator $k$. In-fact, any rescaling $k\to a k$ for some $a\in C^\infty(\cN)$, $a>0$, yields integral curves that re-parametrize the family $\{\beta_p\}_{p\in\Si_0}$. Any change in `base' $\Si_0$ would induce parameter translations for the family $\{\beta_p\}_{p\in\Si_0}$. We will fix our choice of $k$ and for convenience we will also assume $({\La}_-,{\La}_+)\times\mathbb{S}^2\cong\cN$, by discarding all other points.

Throughout the paper we will consider a re-scaling of $k$, which we will again denote $\ubar L := ak$ for some $a>0$, $a\in C^\infty(\cN)$. We conclude therefore that:
\eq{D_{\ubar L}\ubar L = \ka\ubar L,\q \ka = \frac{da}{d\la}.}
Integral curves of $\ubar L$ are denoted by $\ubar\be_p(s)$ for $p\in\Si_0$, whereby $\ubar\be_p(0) = p$ again. By an analogous analysis as above we get an embedding $({\ubar\La}_-,{\ubar\La}_+)\times\Si_{0}\hookrightarrow\cN$.
We denote the null hypersurface associated to the image of $(0,\ubar{\La}_{+})\times \Si_{0}$ by $\Om$ and call it the {\it{null hypersurface generated by $\ubar L$ and $\Si_{0}$}}.
  We also obtain a canonical projection $\pi\cn \Om\to\Si_0$ characterized by the property $\pi(\ubar\beta_p(s)):=p$.  The flow parameter extends to a coordinate function $s\in C^\infty(\Om)$ such that $\partial_s = \ubar L$, $s|_{\Si_0}=0$. The level sets 
\eq{\Si_{s_0}:=\{s=s_0\}\sub\Om} of the function $s$ are diffeomorphic to $\Si_0$ under the projection $\pi$ and we refer to the foliation
\eq{\Om = \bigcup_{s\in(0,{\ubar\La}_+)}\Si_s}
as the {\it{background foliation of $\Om$.}} We denote by $L_s\in\Ga(T^{\perp}\Si_{s})$ the unique null partner of $\ubar L$ satisfying
\eq{\ip{L_s}{L_s}=0,\q\ip{L_s}{\ubar L}=2,\q L_s\perp\Si_s.}
We call $\{\ubar L, L_s\}$ a {\it{null basis}} for $T^{\perp}\Si_{s}$ and also say that $L_s$ is {\it{complementing}} $\ubar L$ in $\Ga(T^\perp\Si_s)$. For geometric quantities along $\{\Si_s\}_s$, we write $\ga_s$ for the induced metric
\eq{\chi_s(V,W) &= -\ip{\mrm{II}_s(W,V)}{L_s} = \ip{D_VL_s}{W}}
for the second fundamental form and 
\eq{\tau_s(V) &= \tfr12\ip{D_V\ubar L}{L_s}}
for the torsion.
By allowing any sections $V,W\in \Ga(T\Om)$ in these formulas, the tensors $\chi_s,\tau_s$ can naturally be extended to $\Om$ and we denote these extensions by $\chi$ and $\tau$ respectively. We also drop the subscript from $\ga_{s}$ and $L_{s}$ and simply write $\ga$ and $L$ instead. The symbol $\n$ denotes the Levi-Civita connection on $\Si_{s}$ and also the gradient of a function with respect to $\ga_{s}$, while $\De$ denotes the Laplace operator. For a tensor $T$ on a Riemannian manifold $(\Si,\si)$ the tensor $\hat T$ denotes the traceless part of $T,$
\eq{\hat T=T-\tfr 12 (\tr_{\si}T)\si.}

Finally, we make a comment on norms of ambient quantities, in particular as arising in the main theorems. As the induced metric $\ga = g_{\cN}$ is degenerate, it does not induce a norm on tensors. However, for tensors annihilated by $\ubar L$, e.g. $\ubar\chi$ and $\ubar\al$, we can define such norms with respect to the background foliation in the following sense:  
\eq{\label{norm}\abs{\hat{\ubar\chi}}^{2}(s,z) = \ga_{s}^{ij}\ga_{s}^{kl}\hat{\ubar\chi}_{ik}\hat{\ubar\chi}_{jl}, }
where we identified $\cN$ as in \eqref{null diffeo} and where coordinates are taken with respect to a local frame  $(e_{i})$ for $\Si_{s}$. A similar definition applies to $\hat{\ubar\al}$. Note that we use the summation convention throughout the paper.

\subsection*{Graphs in null hypersurfaces}
Suppose that $\Si_{\om}\sub\Om$ is the embedding 
of a graphical spacelike surface,
\eq{\label{graph}\Si_{\om}=\{(\om(z),z)\cn z\in \Si_{0}\}}
for some function $\om\in C^{\8}(\Si_{0})$, ${\ubar\La}_-<\om<{\ubar\La}_+$, where the first component in $(\om(z),z)$ denotes the flow parameter of the background foliation. We extend $\om$ constantly along integral curves of $\ubar L$ to a function on $\Om$. The induced geometric quantities of $\Si_{\om}$ are denoted by $L_{\om}$, $\ga_{\om}$, $\s\n$, $\chi_\om$ and $\tau_{\om}$ for the null section complementing $\ubar L$ in $\Ga(T^{\perp}\Si_{\om})$, the  induced metric, its Levi-Civita connection, the second fundamental form with respect to $L_{\om}$, and the torsion with respect to $\{\ubar L,L_{\om}\}$.
Note that due to the annihilating property of $\ubar L$ for $\ubar\chi$ and $\ubar\al$, the definition of the norm \eqref{norm} does not depend on the cross-section, i.e. is the same on any such graph, 
\eq{\abs{\hat{\ubar\chi}}^{2}(\om(z),z) = \ga^{ij}\ga^{kl}\hat{\ubar\chi}_{ik}\hat{\ubar\chi}_{jl}=\ga^{ij}_{\om}\ga^{kl}_{\om}\hat{\ubar\chi}_{ik}\hat{\ubar\chi}_{jl}, }
where we have already used our convention to drop the subscript from $\ga_{s}$ and where indices for the terms on the right hand side are respect to the local frame of $\Si_{\om}$. For quantities defined solely on the graph, such as $\sn\om$ or $\sn^{2}\om$, there is no ambiguity and we use the standard definition of norms, e.g.
\eq{\abs{\sn\om}^{2}= \ga_{\om}^{ij}\del_{i}\om\del_{j}\om.}

\subsection*{Mean curvature flow}
For $\Si_{\om}$ as above the Gaussian formula is
\eq{\label{Gauss}D_{V}W=\s\n_{V}W-\tfr12\chi_\om(V,W)\ubar L-\tfr12\ubar\chi(V,W)L_{\om}\q\fa V,W\in \Ga(T\Si_{\om})}
and hence the correct definition for the mean curvature vector of $\Si_{\om}$ in $\Om$ is
\eq{\cX=\pr_{T\Om}(\vec H)=-\tfr 12(\tr_{\ga_{\om}}\chi_\om) \ubar L.}
For convenience, we simply denote $\tr\chi_\om:=\tr_{\gamma_\om}\chi_\om.$ Consequently, for $T^{*}>0$ the {\it{mean curvature flow}} of spacelike surfaces in $\Om$ is a family 
\eq{x\cn [0,T^{*})\x \bbS^{2}\ra \Om}
of embeddings $x(t,\cdot)$ satisfying
\eq{\label{MCF}\dot x=\cX=-\tfrac12\tr\chi_{\om_t} \ubar L.}
Denote by $(t,\xi)$ elements of $[0,T^{*})\x\bbS^{2}$ and suppose all flow surfaces are given as graphs
\eq{x(t,\bbS^{2})=\Si_{\om_{t}}=\{(\om(t,z(x(t,\xi))),z(x(t,\xi)))\cn \xi\in \bbS^{2}\},}
where $z$ denotes the projection onto $\Si_{0}$, i.e. $z = \pi|_{\Si_{\om_t}}$. Then differentiating
\eq{\om(t,z(x(t,\xi)))=s(x(t,\xi))}
and using $\del_{s} = \ubar L$ gives
\eq{\label{PDE}\fr{\del\om}{\del t}=-\frac12\tr\chi_{\om_t}.}
The right hand side can be expressed in terms of $\om$, see \cite[Lemma~5.1.1]{Roesch:/2017}. Namely we have
\eq{\label{graph-H}\tfr12\tr\chi_{\om_t}=\tfr 12\tr\chi-2\tau(\n\om)-\De\om+2\hat{\ubar\chi}(\n\om,\n\om)+\tfr 12\tr\ubar\chi\abs{\n\om}^{2}-\ka\abs{\n\om}^{2},}
where the right hand side is evaluated at $(\om(t,z),z)$.  This makes \eqref{PDE} a scalar parabolic equation associated to the mean curvature flow, where 
\eq{\label{omega}\om\cn [0,T^{*})\x \Si_{0}\ra ({\ubar\La}_-,{\ubar\La}_+).}

On the other hand if we can prove that for a given function $\om_{0}$, which describes a spacelike surface $x_{0}\cn \bbS^{2}\hra \Si_{\om_{0}}\sub\Om$, we have a maximal solution $\om$ of \eqref{PDE}, then the surfaces $\Si_{t}$ given by
\eq{x(t,\xi)=(\om(t,z(x_{0}(\xi))),z(x_{0}(\xi)))}
solve mean curvature flow. Hence we may assume that $0<T^{*}\leq \8$ is the maximal time of smooth spacelike existence for both equations.

\section{Estimates}\label{Estimates}
To show that the graphical mean curvature flow exists for all times, we need gradient estimates for the function $\om$, and hence we have to differentiate equation \eqref{PDE}. In order to capture the geometric nature of the problem, it is favorable to express $\tr\chi_\om$ in terms of geometric quantities on the hypersurface:

\begin{lemma}\label{l3.1}
For a graphical spacelike hypersurface of $\Om$, there hold
\enum{
\item \eq{\label{L}L_{\om}=L+\abs{\sn\om}^{2}\ubar L-2\sn\om,}
For each $s$, if we identify $\Si_s$ with $\Si_0$ under the induced diffeomorphism $\pi|_{\Si_s}$:
\item {\eq{\label{chi}\chi_\om&=\pi^\ast(\chi)-2(d\om\otimes\pi^\ast(\tau)+\pi^\ast(\tau)\otimes d\om)-2\ka d\om\otimes d\om\\
				&\hp{=}+\abs{\s\nabla\om}^{2}\ubar\chi-2\sn^{2}\om.}}
\item \eq{\tfr12\tr\chi_\om=-\s\De \om-2\tau(\sn\om)+\tfr12\tr\chi+(\tfr 12\tr\ubar\chi+\ka)\abs{\sn\om}^{2}.}
}
\end{lemma}

\pf{
For a tangent vector $V$ of a cross section $\Si_{s}$ at $s=\om(z)$ let 
\eq{\ti V=V+d\om(V)\ubar L,}
which defines an isomorphism of the tangent space of the cross section to that of the graph. 

(i)~The formula for $L_{\om}$ can be checked from the conditions 
\eq{\ip{\ubar L}{L_{\om}}=2,\q \ip{L_{\om}}{L_{\om}}=0,\q \ip{L_{\om}}{\ti V}=0\q\fa \ti V\in \Ga(T\Si_{\om}).}

For (ii) and (iii),
by definition there holds 
\eq{d\om(\ubar L) = 0,\q\text{and}\q [\ubar L,V] = 0.}
Hence
\eq{\chi_\om(\ti V, \ti W)&=\ip{D_{\ti V}(L+\abs{\sn\om}^{2}\ubar L-2\sn \om)}{\ti W}\\
			&=\ip{D_{V+d\om(V)\ubar L}L}{W+d\om(W)\ubar L}\\
			&\hp{=}+\abs{\sn\om}^{2}\ip{D_{\ti V}\ubar L}{\ti W}-2\ip{D_{\ti V}\sn \om}{\ti W}\\
			&=\chi(V,W)-2d\om(V)\tau(W)-2d\om(W)\tau(V)-2\ka d\om(V)d\om(W)\\
			&\hp{=}+\abs{\sn\om}^{2}\ubar\chi(\ti V,\ti W)-2\sn^{2}\om(\ti V,\ti W)\\
			&=\chi(V,W)-2d\om(\ti V)\tau(\ti W)-2d\om(\ti W)\tau(\ti V)+2\ka d\om(\ti V)d\om(\ti W)\\
			&\hp{=}+\abs{\sn\om}^{2}\ubar\chi(\ti V,\ti W)-2\sn^{2}\om(\ti V,\ti W).			}

We observe (ii) from the penultimate equality above. For (iii), note that for an orthonormal frame $(e_{i})$ of a cross section, the frame $(\ti e_{i})$ is also orthonormal. Hence taking the trace of the final equality gives (iii). 
}

We will also need the following famous propagation equations, known as the Raychaudhuri optical equations:

\begin{lemma}\label{l3.02}
For the foliation $\{\Si_s\}_s$ of $\Om$ the following equations hold, where $\pounds$ denotes the Lie derivative:
\eq{\pounds_{\ubar L}\ubar\chi &= -\ubar{\al}+\tfr 12\abs{\hat{\ubar\chi}}^{2}\ga+\tr\ubar\chi \ubar{\hat\chi}+\tfr 14(\tr\ubar\chi)^{2}\ga + \ka\ubar\chi\\
d\tr\ubar\chi(\ubar L) &= -\tfr12(\tr\ubar\chi)^2-\abs{\hat{\ubar\chi}}^2-G(\ubar L,\ubar L)+\ka\tr\ubar\chi,\\
\pounds_{\ubar L}\ubar{\hat\chi} &= -\ubar{\hat\al} +\abs{\hat{\ubar\chi}}^{2}\ga + \ka\hat{\ubar\chi}.} 
\end{lemma}

\pf{For the first two equations, see, for example \cite[Equ.~(4.3), (4.4)]{Roesch:/2017}. For the third:
\eq{\pounds_{\ubar L}\ubar{\hat\chi} &= \pounds_{\ubar L}(\ubar\chi -\tfr 12 \tr{\ubar\chi}\ga)= -\ubar{\hat\al} +\abs{\hat{\ubar\chi}}^{2}\ga + \ka\hat{\ubar\chi}, 				 }
where we have used \cite[Equ.~(4.2)]{Roesch:/2017} and  $\tr{\ubar\al} = G(\ubar L,\ubar L).$
}

Along mean curvature flow, we need several evolution equations.

\begin{lemma}\label{ev-gamma}
Along the mean curvature flow \eqref{MCF} we have the following evolution equations.
\enum{
\item The induced metrics
\eq{\ga_{\om}=x^{\ast}\ip{\cdot}{\cdot}} of the flow hypersurfaces evolve according to
\eq{\pounds_{\del_{t}}\ga_{\om}=-\tr\chi_\om x^{\ast}\ubar\chi.}
\item
The vector $L_{\om}$ evolves according to
\eq{D_{\dot{x}}L_{\om}=\sn \tr\chi_\om+\tr\chi_\om\tau_{\om}^{\sharp}+\tfr12\ka \tr\chi_\om L_{\om}.}
\item
The function $\tr\chi_\om$ evolves according to
\eq{\label{Ev-H}\del_{t}\tr\chi_\om&=\s\De \tr\chi_\om+\tfr12\ka (\tr\chi_\om)^{2}+\tfr 12 \tr\chi_\om\ubar\chi^{i}_{j}(\chi_\om)_{i}^{j}+\tr\chi_\om\s\dive~\tau_{\om}\\
						&\hp{=}+\tfr 12 \tr\chi_\om\ga_{\om}^{ij}\ip{R_{\ubar L x_{i}}L_{\om}}{x_{j}}+2d\tr\chi_\om(\tau_{\om}^{\sharp})+\tr\chi_\om\abs{\tau_{\om}}^{2}.}
}
\end{lemma}

\pf{Fix a local coordinate frame $\del_{i}$ on $\bbS^{2}$ and write
\eq{x_{i}=x_{\ast}(\del_{i}),}
 then 
\eq{\del_{t}\ga_{\om}(\del_{i},\del_{j})&=-\tfr12\ip{D_{x_{i}}(\tr\chi_\om \ubar L)}{x_{j}}-\tfr12\ip{x_{i}}{D_{x_{j}}(\tr\chi_\om \ubar L)}=-\tr\chi_\om \ubar\chi(x_{i},x_{j}).}

To calculate $D_{\dot x}L_{\om}$, note that
\eq{0=\ip{L_{\om}}{L_{\om}}=\ip{D_{\dot x}L_{\om}}{L_{\om}},}
 from $\ip{\ubar L}{L_{\om}}=2$ we get
\eq{0=\ip{D_{\dot x}\ubar L}{L_{\om}}+\ip{\ubar L}{D_{\dot x}L_{\om}}=-\tr\chi_\om\ka+\ip{\ubar L}{D_{\dot x}L_{\om}}}
and finally
\eq{\ip{D_{\dot x}L_{\om}}{x_{i}}=\tfr12\ip{L_{\om}}{D_{x_{i}}(\tr\chi_\om \ubar L)}= d\tr\chi_\om(x_{i})+\tr\chi_\om\tau_{\om}(x_{i}).}

For $\tr\chi_\om$, recall its definition
\eq{\tr\chi_\om=\tr_{\ga_{\om}}\ip{D_{(\cdot)}L_{\om}}{\cdot}=\ga_{\om}^{ij}\ip{D_{x_{i}}L_{\om}}{x_{j}},}
where coordinates to a tensor denote components with respect to the basis $(x_{1},x_{2})$.
First we calculate
\eq{\del_{t}\ip{D_{x_{i}}L_{\om}}{x_{j}}&=\ip{D_{\dot{x}}D_{x_{i}}L_{\om}}{x_{j}}+\ip{D_{x_{i}}L_{\om}}{D_{x_{j}}\dot x}\\
						&=\ip{D_{x_{i}}D_{\dot x}L_{\om}}{x_{j}}-\ip{R_{\dot x x_{i}}L_{\om}}{x_{j}}-\tfr12\ip{D_{x_{i}}L_{\om}}{D_{x_{j}}(\tr\chi_\om \ubar L)}\\
						&=\ip{D_{x_{i}}(\sn \tr\chi_\om+\tr\chi_\om\tau_{\om}^{\sharp}+\tfr12\ka \tr\chi_\om L_{\om})}{x_{j}}\\
						&\hp{=}+\tfr12\tr\chi_\om\ip{R_{\ubar L x_{i}}L_{\om}}{x_{j}}+d\tr\chi_\om(x_{j})\tau_{\om}(x_{i})-\tfr12\tr\chi_\om(\chi_\om)^{k}_{i}\ubar\chi_{kj}\\
						&\hp{=}+\tr\chi_\om\tau_{\om}(x_{i})\tau_{\om}(x_{j})\\
						&=\sn^{2}_{x_{i}x_{j}}\tr\chi_\om+\tr\chi_\om\sn_{x_{i}}\tau_{\om}(x_{j})+\tfr12\ka \tr\chi_\om (\chi_\om)_{ij}\\
						&\hp{=}+\tfr12\tr\chi_\om\ip{R_{\ubar L x_{i}}L_{\om}}{x_{j}}+d\tr\chi_\om(x_{j})\tau_{\om}(x_{i})\\
						&\hp{=}+d\tr\chi_\om(x_{i})\tau_{\om}(x_{j})-\tfr12\tr\chi_\om(\chi_\om)^{k}_{i}\ubar\chi_{kj}+\tr\chi_\om\tau_{\om}(x_{i})\tau_{\om}(x_{j}).}
Hence
\eq{\del_{t}\tr\chi_\om&=\tfr 12 \tr\chi_\om\ubar\chi^{ij}(\chi_\om)_{ij}+\s\De \tr\chi_\om+\tr\chi_\om\s\dive~\tau_{\om}+\tfr12\ka (\tr\chi_\om)^{2}\\
					&\hp{=}+\tfr 12 \tr\chi_\om\ga_{\om}^{ij}\ip{R_{\ubar L x_{i}}L_{\om}}{x_{j}}+2d\tr\chi_\om(\tau_{\om}^{\sharp})+\tr\chi_\om\abs{\tau_{\om}}^{2}.}
}

\begin{cor}
Under the assumptions of \Cref{thm:main}, the mean curvature flow preserves positive $\tr\chi_\om$ up to $T^{*}$.
\end{cor}

\pf{
This follows from the strong maximum principle applied to the evolution equation 
\eqref{Ev-H}.
}

Now the goal is to deduce an estimate for the norm of the gradient
\eq{u:=\tfr 12\abs{\sn\om}^{2}.}
We first establish the $C^{0}$-estimates.

\begin{prop}\label{comp}
During the evolution the flow ranges in the fixed compact domain enclosed by $\Si_{0}$ and $\Si_{\om_{0}}$.
\end{prop}

\pf{
This follows from the maximum principle applied to the evolution of the graph function
\eq{\del_{t}\om=-\tfr12\tr\chi_{\om_{t}}=\s\De \om+2\tau(\sn\om)-\tfr12\tr\chi-(\tfr 12\tr\ubar\chi+\ka)\abs{\sn\om}^{2},}
since $\Si_{0}$ has $\tr\chi\leq 0$ and $\Si_{\om_{t}}$ has $\tr\chi_{\om_t}>0$.
}
We define the parabolic operator
\eq{\cP=\del_{t}-\s\De-2\tau(\sn(\cdot)).}
In the following, $\cO(y)$ denotes any function that is bounded up to order one, when $y\ra \8$, i.e. the estimate
\eq{\limsup_{\abs{y}\ra \8}\fr{\abs{\cO(y)}}{\abs{y}}\leq C,}
where $C$ depends on the data of the problem (i.e. on $\Om$ and $\Si_{0}$).

\begin{lemma}\label{ev-grad}
Along \eqref{PDE} the quantity
\eq{u:=\tfr 12\abs{\sn\om}^{2}}
satisfies the equation
\eq{\cP u &= \cO(\abs{\sn\om}^3)+\tfr12\tr\chi_\om\hat{\ubar\chi}(\sn\om,\sn\om)-\abs{\sn^2\om}^2-2u\tr(\hat{\ubar\chi}\circ\sn^2\om)\\
			&\hp{=}-\tr\ubar\chi du(\sn\om)+\Big((\tr\ubar\chi)^2+4\abs{\hat{\ubar\chi}}^2-2\ka\tr\ubar\chi+4d\ka(\ubar L)\Big)u^{2}\\
			&\hp{=}-2\ka\hat{\ubar\chi}(\sn\om,\sn\om)u+2\hat{\ubar\al}(\sn\om,\sn\om)u.}
\end{lemma}

\pf{
We have to differentiate the PDE
\eq{\label{ev-grad-1}\del_{t}\om=-\tfr12\tr\chi_\om=\s\De \om+2\tau(\sn\om)-\tfr12\tr\chi-(\tr\ubar\chi+2\ka)u}
covariantly with respect to $\ga_{\om} = x^{\ast}g$ in the variable $z\in \Si_{0}$. Pick a local coordinate frame $(\del_{i})_{1\leq i\leq 2}$ for $\Si_{0}$ and denote by $(x_{i})$ the induced frame on the graph given by
\eq{x_{i}=\del_{i}+\del_{i}\om \ubar L.}
We evaluate the terms in the following expression separately:
\eq{-\tfr12d\tr\chi_\om(\sn\om)&=d(\s\De \om+2\tau(\sn\om)-\tfr12\tr\chi-(\tr\ubar\chi+2\ka)u)(\sn\om).}

(i)~Firstly we have
\eq{&-\tfr12d\tr\chi(\sn\om)-d((\tr\ubar\chi+2\ka) u)(\sn\om)\\
=~&-\tfr12d\tr\chi(\sn\om)-d\tr\ubar\chi(\sn\om)u-2 d\ka(\sn\om)u-(\tr\ubar\chi+2\ka) du(\sn\om)\\
=~&\cO(\abs{\sn\om}^{3})+2(\tfr12(\tr\ubar\chi)^2+\abs{\hat{\ubar\chi}}^2+G(\ubar L,\ubar L)-\ka\tr\ubar\chi)u^2\\
&-4d\ka(\ubar L)u^2 - (\tr\ubar\chi+2\ka)du(\sn\om)\\
=~&\cO(\abs{\sn\om}^3)+\Big((\tr\ubar\chi)^2+2\abs{\hat{\ubar\chi}}^2+2G(\ubar L,\ubar L)-2\ka\tr\ubar\chi-4d\ka(\ubar L)\Big)u^{2}\\
& - (\tr\ubar\chi+2\ka)du(\sn\om),}
having used \Cref{l3.02} and
\eq{\label{slashnabla}\sn\om = \n\om + \abs{\sn\om}^{2}\ubar L = \n\om+2u\ubar L}
 in the penultimate equality.

(ii)~Secondly there is
\eq{d(\s\De\om)(\sn\om)=\s\De u-\abs{\sn^{2} \om}^{2}-\s\Rc(\sn \om,\sn\om).}
We notice that we deal with surfaces and hence the Ricci curvature is 
\eq{\s\Rc=\s{\mathcal{K}}\ga_{\om},}
where $\s{\mc{K}}$ is the Gauss curvature of the graph. Using \cite[Prop.~3.0.1]{Roesch:/2017},
\eq{\label{ev-grad-2}2\s{\mc{K}}&=\tfr12\tr\chi_\om\tr\ubar\chi-\tr(\hat{\ubar\chi}\circ\hat\chi_\om)-S-2G(\ubar L,L_{\om})-\tfr 12\ip{R_{\ubar LL_{\om}}\ubar L}{L_{\om}}.
		}
 We conclude with the help of \eqref{L}, \eqref{chi}, the symmetries of the Riemann tensor, \eqref{slashnabla}
and
\eq{\tr\ubar\al = G(\ubar L,\ubar L),}
that
\eq{d(\s\De\om)(\sn\om)&=\s\De u-\abs{\sn^{2} \om}^{2}-2\s{\mc{K}}u\\
				&=\s\De u-\abs{\sn^{2} \om}^{2}-\tfr12\tr\chi_\om \tr\ubar\chi u+\tr(\hat{\ubar\chi}\circ\hat\chi_\om)u+Su\\
				&\hp{=}+2G(\ubar L,L+\abs{\sn\om}^{2}\ubar L-2\sn\om)u+\tfr 12\ip{R_{\ubar L(L-2\sn\om)}\ubar L}{L-2\sn\om}u\\
				&=\s\De u-\abs{\sn^{2}\om}^{2}-\tfr12\tr\chi_\om\tr\ubar\chi u -2\ka\hat{\ubar\chi}(\sn\om,\sn\om)u+2u^2\abs{\hat{\ubar\chi}}^2\\
				&\hp{=}-2u\tr(\hat{\ubar\chi}\circ\sn^2\om)+\cO(\abs{\sn\om}^{3})-2G(\ubar L,\ubar L)u^{2}+2\hat{\ubar\al}(\sn\om,\sn\om)u.}

(iii)~Finally we have, locally extending $\tau$ to $T^{0,1}(M)$ and using
\eq{\ubar\chi(\sn\om,\sn\om) = \ubar\chi(\n\om,\n\om) = \cO(\abs{\sn\om}^{2}),}
\eq{2d(\tau(\sn\om))(\sn\om)&=2\sn\om(\tau(\sn\om))\\
						&=2D\tau(\sn\om,\sn\om) + 2\tau(D_{\sn\om}\sn\om)\\
						&=\cO(\abs{\sn\om}^3)+8D\tau(\ubar L,\ubar L)u^{2}+ 2\tau(\sn u)\\
						&\hp{=}-\chi_\om(\sn\om,\sn\om)\tau(\ubar L)+ 2\ubar\chi(\sn\om,\sn\om)\tau(\ubar L)u.
						}
Now we use
\eq{D\tau(\ubar L,\ubar L)=\ubar L(\tau(\ubar L)) - \tau(D_{\ubar L}\ubar L)=\ubar L\ka - \ka^{2}}
and
\eq{\chi_\om(\sn\om,\sn\om) &= \ip{D_{\sn\om}(L+2u\ubar L-2\sn\om)}{\sn\om}\\
					&=\ip{D_{\sn\om}L}{\sn\om} + 2\ip{D_{\sn\om}\ubar L}{\sn\om}u - 2\ip{\sn_{\sn\om}\sn\om}{\sn\om}\\
					&=\cO(\abs{\sn\om}^{3})+4\ip{D_{\ubar L}L}{\ubar L}u^{2} + 2\ubar\chi(\sn\om,\sn\om)u - 2\ip{\sn u}{\sn\om}\\
					&=\cO(\abs{\sn\om}^{3})-8\ka u^{2} + 2\ubar\chi(\sn\om,\sn\om)u - 2\ip{\sn u}{\sn\om}}
to deduce
\eq{2d(\tau(\sn\om))(\sn\om)&=\cO(\abs{\sn\om}^{3})+8d\ka(\ubar L)u^{2}+ 2\tau(\sn u) + 2\ka \ip{\sn u}{\sn\om}.
						}

The previous three steps provide all the terms coming from the right hand side of \eqref{ev-grad-1}. For the left hand side we denote by $(\ga_{\om}^{ij})$ the inverse of $\ip{x_{i}}{x_{j}}$ and differentiate $u$ in time:
\eq{\del_{t}u=\tfr 12\del_{t}(\ga_{\om}^{ij}\del_{i}\om\del_{j}\om)&=\tfr12\tr\chi_\om\ga_{\om}^{ik}\ubar\chi(x_{k},x_{l})\ga_{\om}^{lj}\del_{i}\om\del_{j}\om-\tfr12d\tr\chi_\om(\sn\om)\\
					&=\tfr12\tr\chi_\om\ubar\chi(\sn\om,\sn\om)-\tfr12d\tr\chi_\om(\sn\om)\\
					&=\tfr12\tr\chi_\om\hat{\ubar\chi}(\sn\om,\sn\om)+\tfr12\tr\chi_\om\tr\ubar\chi u-\tfr12d\tr\chi_\om(\sn\om).}
					
Collecting the terms from items (i)-(iii), we get
\eq{\del_{t}u	&= \tfr12\tr\chi_\om\hat{\ubar\chi}(\sn\om,\sn\om)+\tfr12\tr\chi_\om\tr\ubar\chi u\\
			&\hp{=}+\cO(\abs{\sn\om}^3)+\Big((\tr\ubar\chi)^2+2\abs{\hat{\ubar\chi}}^2+2G(\ubar L,\ubar L)-2\ka\tr\ubar\chi-4d\ka(\ubar L)\Big)u^{2} \\
			&\hp{=}- (\tr\ubar\chi+2\ka)du(\sn\om)+\cO(\abs{\sn\om}^{3})+\s\De u-\abs{\sn^{2}\om}^{2}-\tfr12\tr\chi_\om\tr\ubar\chi u\\
			&\hp{=} -2\ka\hat{\ubar\chi}(\sn\om,\sn\om)u+2u^2\abs{\hat{\ubar\chi}}^2-2u\tr(\hat{\ubar\chi}\circ\sn^2\om)-2G(\ubar L,\ubar L)u^{2}\\
			&\hp{=}+2\hat{\ubar\al}(\sn\om,\sn\om)u+ \cO(\abs{\sn\om}^{3})+8d\ka(\ubar L)u^{2}+ 2\tau(\sn u) + 2\ka \ip{\sn u}{\sn\om}\\
			&=\cO(\abs{\sn\om}^3)+\tfr12\tr\chi_\om\hat{\ubar\chi}(\sn\om,\sn\om)-\abs{\sn^2\om}^2-2u\tr(\hat{\ubar\chi}\circ\sn^2\om)+\s\De u\\
			&\hp{=}-\tr\ubar\chi du(\sn\om) +2\tau(\sn u)+\Big((\tr\ubar\chi)^2+4\abs{\hat{\ubar\chi}}^2-2\ka\tr\ubar\chi+4d\ka(\ubar L)\Big)u^{2}\\
			&\hp{=}-2\ka\hat{\ubar\chi}(\sn\om,\sn\om)u+2\hat{\ubar\al}(\sn\om,\sn\om)u.}
giving the result.				
}

Now we can turn to the $C^{1}$-estimates.

\begin{prop}\label{prop:grad-bound}
Under the assumptions of \Cref{thm:main}, along the mean curvature flow starting from any embedded spacelike surface in the null hypersurface $\Om$, the function
\eq{u=\tfr 12\abs{\sn\om}^{2}}
is uniformly bounded on $[0,T^{*})\x \Si_{0}$.
\end{prop}

\pf{
We have to control the right hand side in \Cref{ev-grad}. First we deal with the term involving $\tr \chi_{\om}$. We have
\eq{&\tfr12\tr\chi_\om\hat{\ubar\chi}(\sn\om,\sn\om)-\abs{\sn^{2}\om}^{2}-2u\tr(\hat{\ubar\chi}\circ\sn^2\om)\\
	=~&(-\s\De \om-2\tau(\sn\om)+\tfr12\tr\chi+(\tfr 12\tr\ubar\chi+\ka)\abs{\sn\om}^{2})\hat{\ubar\chi}(\sn\om,\sn\om)-\abs{\sn^{2}\om}^{2}\\
	&-2u\tr(\hat{\ubar\chi}\circ\sn^2\om)\\
	=~&\cO(\abs{\sn\om}^3)-\s\De\om\hat{\ubar\chi}(\sn\om,\sn\om)-\abs{\sn^{2}\om}^{2}-2u\tr(\hat{\ubar\chi}\circ\sn^2\om)\\
	&+(\tr\ubar\chi-2\ka) \hat{\ubar\chi}(\sn\om,\sn\om)u.}
	
Since 
\eq{\abs{\s\nabla^2\om}^2 = \abs{\hat{\s\nabla}^2\om}^2+\tfr12(\s\De\om)^2,\q  \tr(\hat{\ubar\chi}\circ\s\nabla^2\om) = \tr(\hat{\ubar\chi}\circ\hat{\s\nabla}^2\om),} we may complete the square on the second, third and forth terms above to conclude
	\eq{&\tfr12\tr\chi_\om\hat{\ubar\chi}(\sn\om,\sn\om)-\abs{\sn^{2}\om}^{2}-2u\tr(\hat{\ubar\chi}\circ\sn^2\om)\\
	\leq~&\cO(\abs{\sn\om}^{3})+\tfr 12\hat{\ubar\chi}(\sn\om,\sn\om)^{2}+u^2\abs{\hat{\ubar\chi}}^2+(\tr\ubar\chi-2\ka)\hat{\ubar\chi}(\sn\om,\sn\om)u\\
	\leq~& \cO(\abs{\sn\om}^3)+3\abs{\hat{\ubar\chi}}^2u^2+(\tr\ubar\chi-2\ka)\hat{\ubar\chi}(\sn\om,\sn\om)u,}
where we have used Cauchy-Schwarz in the final estimate. Hence, under the assumptions of \Cref{thm:main}, and using \Cref{l3.02},

\eq{\cP u &\leq \cO(\abs{\sn\om}^3)-\tr\ubar\chi du(\sn\om)+\Big((\tr\ubar\chi)^2+7\abs{\hat{\ubar\chi}}^2-2\ka\tr\ubar\chi+4d\ka(\ubar L)\Big)u^2\\
&\hp{=}+(\tr\ubar\chi-4\ka)\hat{\ubar\chi}(\sn\om,\sn\om)u+2\hat{\ubar\al}(\sn\om,\sn\om)u\\
&\leq \cO(\abs{\sn\om}^3)-\tr\ubar\chi du(\sn\om)+\Big(5\abs{\hat{\ubar\chi}}^2-2G(\ubar L,\ubar L)+2d(2\ka-\tr\ubar\chi)(\ubar L)\Big)u^2\\
&\hp{=}+(\tr\ubar\chi-4\ka)\hat{\ubar\chi}(\sn\om,\sn\om)u+2\hat{\ubar\al}(\sn\om,\sn\om)u\\
&\leq \cO(\abs{\sn\om}^3)-\tr\ubar\chi du(\sn\om)+((\tr\ubar\chi-4\ka)\hat{\ubar\chi}(\sn\om,\sn\om)\\
&-2\abs{(\tr\ubar\chi-4\ka)\hat{\ubar\chi}}u)u+2(\hat{\ubar\al}(\sn\om,\sn\om)-2\abs{\hat{\ubar\al}}u)u\\
&\leq \cO(\abs{\sn\om}^3)-\tr\ubar\chi du(\sn\om).
}

In order to estimate $u$, we combine this evolution equation with that of $\om$, which is
\eq{\cP\om = -\tfr12\tr\chi - (\tr\ubar\chi+2\ka)u.}
At points where $u>0$, define the auxiliary function 
\eq{\phi = \log u + f(\om),}
with $f$ yet to be determined. Then
\eq{\cP\phi = \fr{\cP u}{u} +\abs{\sn\log u}^{2}+ f'\cP\om - f''\abs{\sn\om}^{2}.}
At maximal points of $\phi$ we have
\eq{\sn\log u =\tfr{1}{u} \sn u= -f'\sn\om}
and hence at such points
\eq{0&\leq \cP\phi\\
	&\leq \cO(\abs{\sn\om})+2\tr\ubar\chi f'u+2f'^{2}u-\tfr12f'\tr\chi-(\tr\ubar\chi+2\ka)f'u-2f''u\\
	&\leq \cO(\abs{\sn\om})-\tfr12f'\tr\chi+\Big(2f'^{2}+(\tr\ubar\chi-2\ka)f'-2f''\Big)u.}

Denoting by $\cC\sub\cN$ the compact region of \Cref{comp}, we set
\eq{\la = \max\{1,\sup_{\cC}(\tr\ubar\chi-2\ka)\}}
and with $\ep = \tfr 12 e^{-\la\max \om_{0}}
$ define
\eq{
f(\om) = -\log(e^{-\la \om}-\ep) - \la\om.
}
Then
\eq{f' = \fr{\la e^{-\la\om}}{e^{-\la\om}-\ep} -\la,\q f'' &= -\fr{\la^{2}e^{-\la\om}}{e^{-\la\om}-\ep} + \fr{\la^{2}e^{-2\la\om}}{(e^{-\la\om}-\ep)^{2}}}
and
\eq{f''&= -\la f'-\la^{2} + f'^{2} + 2\la f'+\la^{2} = f'^{2} + \la f'.  }
Hence
\eq{0&\leq \cO(\abs{\sn\om})-\tfr12f'\tr\chi +\Big(\tr\ubar\chi-2\ka-2\la \Big)f'u\\
&\leq \cO(\abs{\sn\om})-\la f'u-\tfr12f'\tr\chi.}
This estimate gives a contradiction, if $u$ is too large. Hence $u$ is uniformly bounded in terms of the data of the problem.	
}

From quasi-linear regularity theory this bootstraps to estimates in $C^{\8}$.

\begin{prop}\label{Ck}
Under the assumptions of \Cref{thm:main}, the mean curvature flow starting from any embedded spacelike outer un-trapped surface $\Si_{\om_{0}}$ exists for all times and satisfies uniform estimates in any $C^{k}(\Si_{0})$ norm, where $\Si_{0}$ is understood to be equipped with its induced metric.
\end{prop}

\pf{
Due to the assumption, the functions $\om(t,\cdot)$ are uniformly bounded independently of $t<T^{*}$. Since 
\eq{\abs{\n\om}=\abs{\sn\om},}
the differential $d\om$ is bounded with respect to the induced metric of $\Si_{s=\om}$, which is in turn equivalent to the induced metric of $\Si_{0}$, since the surfaces $\Si_{\om}$ range in a compact set of $\Om$.
Hence $\om$ solves the quasi-linear equation 
\eq{\fr{\del\om}{\del t}(t,z)=-\tfr12\tr\chi_\om(t,z),}
to which parabolic regularity theory may be applied, \cite[Ch.~XII]{Lieberman:/1998}, to obtain $C^{1,\al}$ estimates for $\om$. This is sufficient to start a bootstrap argument involving parabolic Schauder estimates to obtain the $C^{k}$-estimates, see \cite{Lieberman:/1998}, up to time $T^{*}$. The standard result on short-time existence, see for example \cite[Thm.~2.5.7]{Gerhardt:/2006}, implies that the flow can be continued beyond any finite time and hence we obtain $T^{*}=\8$.
}

\subsection{Completion of the proof of \texorpdfstring{\Cref{thm:main}}{Theorem 1.1}}\label{proof main}

The mean curvature flow starting from $\Si_{\om_{0}}$ exists for all times, and the graph functions of its leaves $\Si_{\om_{t}}$ are strictly decreasing due to 
\eq{\fr{\del\om}{\del t}(t,z)=-\tfr12\tr\chi_{\om_t}}
and the preservation of $\tr\chi_{\om_t}>0$. Hence the pointwise limit
\eq{\om_{\8}(z)=\lim_{t\ra\8}\om(t,z)}
exists and the convergence takes place in the topology of $C^{\8}$, due to the $C^{k}$-estimates from \Cref{Ck}. Integration gives
\eq{\8>\abs{\om_{\8}(z)-\om_{0}(z)}=\tfr12\int_{0}^{\8}\tr\chi_{\om}(t,z)~dt}
and hence, due to the regularity estimates, the limit is a MOTS:
\eq{\lim_{t\ra \8}\tr\chi_{\om}(t,\cdot)=0.}

\section{Null Cones}\label{sec:cone}
In this section we wish to prove \Cref{thm:cone} as our first application of \Cref{thm:main} under the energy condition (\ref{energy}). First we need to specialize to a specific null geometry:
\begin{defn}
	Let $\cN$ be a null hypersurface as described in \Cref{sec:setup}. If the second fundamental form of $\cN$ associated to the vector field $k$ satisfies the condition $\tr K>0,$
we say $\cN$ is a Null Cone.	 
\end{defn}
In general relativity the Einstein tensor $G$ is coupled to the stress-energy tensor of matter as modelled by a spacetime $(M,g)$. A somewhat weak consequence of the physical assumption of a non-negative energy density distribution within $(M,g)$ is the so-called \textit{null convergence condition}
\eq{G(k,k)\geq 0.}
This condition gives rise to Null Cone structures in $(M,g)$. Under fairly generic assumptions on $\cN$, we support this claim with the following well known result:

\begin{lemma}
 Suppose $\La_+=\infty$, specifically $\cN\cong (\La_-,\infty)\times\bbS^2$, and satisfies the null convergence condition. Then $\tr K\geq 0$. If, in addition, the set $E:=\{p\in\cN|\tr K(p)>0\}$ admits a surjection $\pi\cn E\to \Si_0$, then $\cN=E$. 
\end{lemma}
\pf{

If we assume for a contradiction that $\tr K(p)<0$ for some $p\in\cN$, then taking the geodesic $\be_p(\la)$ associated to $k$, and applying \Cref{l3.02} to the case $\ubar L = k$ (giving $\ka\equiv 0$) we observe
\eq{\frac{d}{d\la}\frac{1}{\tr K} = \fr12+\fr{\abs{\hat{K}}^2+G(k,k)}{(\tr K)^2}\geq \fr12,}
as long as $\tr K(\la)<0$. We conclude $\tfr{1}{\tr K}(\la)\geq\tfr12 \la+\tfr{1}{\tr K}(p)$, so that $\tfr{1}{\tr K}(\la)\to 0$ as $\la\to\la_\star^-$ for some $\la_\star>0$. Consequently, $\tr K(\la)\to -\infty$ as $\la\to \la_\star^-$ contradicting the smoothness assumption of $\cN$ under $\La_+ = \infty$. It follows that $\tr K\geq 0$ throughout $\cN$. If we take $p\in E$, we have for $\la\geq 0$
\eq{\fr{1}{\tr K}(\la) \geq \tfr12\la+\fr{1}{\tr K}(p)>0,}
which implies $\tr K(\la)>0$. Alternatively, for $\la\leq 0$
\eq{\infty> \fr{1}{\tr K}(p) \geq -\tfr12\la+\fr{1}{\tr K}(\la),}
again implies $\tr K(\la)>0$. Since every $q\in\cN$ admits a unique geodesic along $k$ intersecting $\Si_0$, we conclude the same geodesic intersects the set $E$ under our second hypothesis giving $\tr K(q)>0$.
}

A consequence of the above result is that any 2-sphere $\Si\sub M$ satisfying the outer un-trapped condition to the timelike past
\eq{-\ip{k}{\vec{H}}=\tr K>0,}
for $k\in \Ga(T^\perp\Si)$ a past-pointing null vector field, forms a Null Cone geometry if the congruence of null geodesics along $k$ extend infinitely. Moreover, by the \textit{first variation of area formula}, the variation of the area form $dA$ on $\Si$ is
\eq{\de_k(dA) = \tr KdA.}
It follows that a Null Cone is foliated by 2-spheres with pointwise expanding area form along $k$, hence the name \textit{Null Cone}. In physical spacetimes we observe that any infinite and \textit{shear-free} Null Cone, i.e. satisfying $\hat{K}\equiv 0$, satisfies (\ref{energy}) as a consequence of the null convergence condition. We note that in the shear-free case also $\hat\al$ vanishes, see \Cref{l3.02}. Consequently, the mean curvature flow can be initiated from any outer un-trapped surface irrespective of it's proximity to $\Si_0$ within a shear-free Null Cone.

\subsection{Proof of \texorpdfstring{\Cref{thm:cone}}{Theorem 1.2}}
Our goal is to show the existence of a vector field $\ubar L = ak$ satisfying the gauge condition of \Cref{thm:main}. This translates into showing the existence of $a\in C^\infty(\cC)$, $a>0$, for the region $\cC$ bounded by $\Si_0$, and $\Si_{\om_{0}}$ satisfying the gauge conditions of \Cref{thm:main}, and such that the integral curves along $\ubar L$ cover all of $\cC$. For convenience, we will denote $\be:=\tfr12\tr K$, and for a function $f$, $f':=df(k)$.

We start, for some $v_0\in C^\infty(\Si_0)$, $0<v_0<1$, by solving the ODEs
\eq{v' = \be(1-v)^2,\q v(0,p) = v_0(p).}
From standard existence results, and a phase space analysis we observe:
\eq{0<v(\la,p)<1,\q \la\in[0,\La_+),}
in-fact,
\eq{v(\la,p)=\fr{v_0(p)+(1-v_0(p))\int_0^\la\beta(u,p)du}{1+(1-v_0(p))\int_0^\la\beta(u,p)du}.}
Combining $0<v<1$ with (\ref{energy}) and then using \Cref{l3.02} we have
\eq{\frac{v'}{\be}&\geq (1-v)^2-\fr{v}{2\be^2}(\abs{\hat{K}}^2+G(k,k)-4\be\abs{\hat{K}})+\fr{v^2}{2\be^2}(\tfr{7}{2}\abs{\hat{K}}^2+2\abs{\hat{\al}})\\
&=1-v+v^2-v\br{\tfr{1}{\be}}'+2\abs{\hat{K}}\fr{v}{\be}+\fr{v^{2}}{\be^{2}}(\tfr74\abs{\hat{K}}^2+\abs{\hat{\al}}),}
giving
\eq{\fr{\be^{2}}{v^{2}}\br{\br{\fr{v}{\be}}'-v^{2}+v-1} &\geq \tfr74\abs{\hat{K}}^2+\abs{\hat{\al}}+\be\abs{\hat{K}}+(2\tfr{\be}{v}-\be)\abs{\hat{K}}\\
&\geq \tfr74\abs{\hat{K}}^2+\abs{\hat{\al}}+\abs{2\tfr{\be}{v}-\be}\abs{\hat{K}}.}
We claim that our desired $a(\la,p)$ is given by solving the family of ODEs
\eq{(\log a)' = \be v^{-1},\q a(0,p) = 1.}
The differential inequality above combined with \Cref{l3.02} then gives
\eq{&\tfr12G(k,k)-\fr{d}{d\la}((\log a)'-\be)-(\log a)'((\log a)'-\be )\\
	\geq~& \abs{2(\log a)'-\be}\abs{\hat{K}}+\fr54\abs{\hat{K}}^2+\abs{\hat{\al}}.}
Multiplying throughout by $2a^2$ and simplifying, we finally observe
\eq{G(\ubar L,\ubar L)-a(2a'-\tr\ubar\chi)'\geq \abs{4a'-\tr\ubar\chi}\abs{\hat{\ubar\chi}}+\fr52\abs{\hat{\ubar\chi}}^2+2\abs{\hat{\ubar\al}},}
since $\ka = a'$, and $ak = \ubar L$, we've satisfied the gauge condition of \Cref{thm:main}. We also observe from the fact that $\infty>a(\la,p)\geq 1$ throughout $\cC$, that $\ubar L$ is globally defined throughout.

\section{Existence of a global outer un-trapped foliation}\label{sec:foliation}
In physical applications, where $\Omega$ admits an infinite \textit{asymptotically flat} region, for example the study of the Penrose inequality \cite{Penrose:/1973} or mass in spacetimes \cite{ChruscielPaetz:04/2014},  one observes a foliation of $\Om$ by future un-trapped surfaces in a neighborhood of infinity. In-fact, such foliations are plentiful in a neighborhood of infinity (see for example, \cite{MarsSoria:09/2015}). For $\Om$ also admitting an outermost MOTS $\Si_{out}$, it is physically reasonable to also expect this future un-trapped foliation to extend all the way up-to $\Si_{out}\sub \Om$ (see, for example \cite{Sauter:/2008}). More specifically, if $\Sigma_{out} := \{s=\omega(0,z)\}$, and $\sigma> 0$, we expect:
\eq{\frac{\partial}{\partial\sigma}\omega(\sigma,z) > 0,\q \tr\chi_{\om}(\sigma,z) > 0.}
In this section, we apply \Cref{thm:main} to the assumption of an un-trapped foliation in a neighborhood of $\La_+$ (note here $\La_+ = \infty$ is a possibility), to show the existence of a global future un-trapped foliation to the timelike past of the MOTS of convergence. Moreover, also denoted $\Si_{out}\sub\Om$, we show that this MOTS is indeed outermost. 

Under the existence of a future un-trapped cross-section, we reparametrize our global coordinate chart on $\Om$ so that a leaf $\Si_{\La}$, for some $\La\in(0,\La_+)$, corresponds with this cross-section. We now also assume that all $\{\Si_s\}_{\La\leq s<\La_+}$ are future un-trapped, in addition to the hypotheses of \Cref{thm:main}. Thus we may initiate mean curvature flow at $\ti\omega(0,z) = \Lambda$ obtaining the smooth solution $\ti\omega(t,z)$ for $t>0$. Here $\ti\om$ arises from the solution $\om$ of \eqref{PDE}, which arose from fixing the gauge $\ubar L$ according to the hypotheses of \Cref{thm:main}. 
It suffices therefore to show that the two families $\{\Si_s\}_{s> \La}$ and $\{\Si_{\ti\omega_t}\}_{t>0}$ can be `glued and smoothed' at the boundary $\Sigma_\Lambda$ to form a global future un-trapped foliation. We achieve this is by smoothing the continuous function
\eq{v(\lambda,z):=\begin{cases}
         \ti\omega(\Lambda-\lambda,z), & \text{for}~ \lambda<\Lambda\\
         \lambda, & \text{for}~ \lambda\geq \Lambda.
        \end{cases}}
In order to do so, we adapt a standard mollifying procedure that we will now outline. With $C$ chosen so that $\int\eta(s)ds = 1$, we start with the mollifier
\eq{\eta(s):=\begin{cases} 
          Ce^{\frac{1}{s^2-1}}, & \text{for}~ |s|<1 \\
          0, & \text{for}~ |s|\geq 1
          \end{cases}}
and define $\eta_\epsilon(s):=\tfrac{1}{\epsilon}\eta\br{\tfrac{s}{\epsilon}}.$ 
We mollify $v(\lambda,z)$:
\eq{v_\epsilon(\lambda,z): = \int\eta_\epsilon(\lambda-u)v(u,z)du.}
It follows that
\begin{align}
\partial_\lambda v_\epsilon(\lambda,z) &= \int \partial_\lambda\eta_\epsilon(\lambda-u)v(u,z)du,\\
\partial^\alpha_{z}v_{\epsilon}(\lambda,z)&=\int\eta_\epsilon(\lambda-u)\partial^\alpha_{z} v(u,z)du
\end{align}
for any multi-index $\alpha$, and both functions above are smooth. Since
\eq{\partial^\alpha_z v(\lambda,z) = \begin{cases}
\partial_z^\alpha \ti\omega(\Lambda-\lambda,z), &\text{for}~\lambda<\Lambda\\
0,&\text{for}~\lambda\geq \Lambda
\end{cases}} 
is continuous, standard results imply $\partial_z^\alpha v_\epsilon\to \partial_z^\alpha v$ uniformly as $\epsilon\to 0$ on $\mathcal{U}:=(\Lambda-\delta,\Lambda+\delta)\times\mathbb{S}^2$ for some $\delta>0$. From our constraints on $\Omega$, we observe 
\eq{\tr\chi_{v}(\lambda,z)>4c>0} for some $c>0$ throughout $\bar{\mathcal{U}}$. Define the partition of unity $\{\zeta_i\}_{1\leq i\leq 3}$ for $\bbR$,
\eq{\zeta_1(\la)&= \begin{cases}1, & \lambda\leq \Lambda-\delta\\
				0, & \lambda\geq \Lambda-\frac12\delta,
				\end{cases}\\
\zeta_2(\la)&= \begin{cases} 1, & |\lambda-\Lambda|\leq\frac12\delta\\
					0, & |\lambda-\Lambda|\geq \delta
					\end{cases}\\
		\zeta_3(\lambda)&= \zeta_1(2\Lambda-\lambda).}
 We may therefore construct the smooth function
\eq{\ti \omega_\epsilon(\lambda,z):=\zeta_1(\lambda)\ti \omega(\Lambda-\lambda,z)+\zeta_2(\lambda)v_{\epsilon}(\lambda,z)+\zeta_3(\lambda)\lambda}
and we observe that $\partial_z^\alpha \ti\omega_\epsilon\to \partial_z^\alpha v$ uniformly on $\bar{\mathcal{U}}$ as $\epsilon\to 0$. From \eqref{chi}, we observe that the future null expansion $\tr\chi_{\ti\om}$ depends only upon $\partial_z^\alpha \ti \omega_\epsilon$ for $|\alpha|\leq 2$. Consequently, 
\eq{\tr\chi_{\ti\om_{\ep}}(\lambda,z)\to  \tr\chi_{v}(\lambda,z)} uniformly as $\epsilon\to 0$ throughout $\bar{\mathcal{U}}$. We therefore conclude that $\epsilon_0>0$ exists such that $ \tr\chi_{\ti\om_{\ep}}(\lambda,z)>c$ throughout $\bar{\mathcal{U}}$ for $\epsilon\leq \epsilon_0$,
 \eq{ \tr\chi_{\ti\omega_{\epsilon}}(\lambda,z) =  \tr\chi_{\ti\omega}(\Lambda-\lambda,z)>0} on $\lambda\leq \Lambda-\delta$ and
 \eq{ \tr\chi_{\ti\omega_{\epsilon}}(\lambda,z) =  \tr\chi|_{\Sigma_\lambda}>0} on $\lambda\geq \Lambda+\delta$. Finally, we wish to show $\partial_\lambda\ti\omega_\epsilon>0$ for sufficiently small $\epsilon$. From the fact that \eq{\ti\omega_\epsilon(\lambda,z)\equiv \ti\omega(\Lambda-\lambda,z)} for $\lambda\leq \Lambda-\delta$, and $\ti\omega_\epsilon\equiv \lambda$ for $\lambda\geq \Lambda+\delta$ it will suffice to show $\partial_\lambda\ti\omega_\epsilon>0$ on $\bar{\mathcal{U}}$ for sufficiently small $\epsilon$. From a simple integration by parts:
\eq{\partial_\lambda v_\epsilon(\lambda,z) = \int_{-\infty}^\Lambda\tfr12\eta_\epsilon(\lambda-u) \tr\chi_{\ti\omega}(\Lambda-u,z)du+\int_\Lambda^{\infty}\eta_\epsilon(\lambda-u)du>0,}
leaving us only the sets 
\eq{\mathcal{U}_{+}:=(\Lambda+\tfr12\delta,\Lambda+\delta)\times\mathbb{S}^2,\,\,\,\cU_-:=(\La-\de,\La-\tfr12\de)\times\mathbb{S}^2} to consider. On $\mathcal{U}_-$, we know $\zeta_1'+\zeta_2'\equiv 0$, so taking $\epsilon<\frac12\delta$ we observe: 
\eq{
\partial_\lambda \ti\omega_\epsilon(\lambda,z) &= \zeta_1'(\lambda)\big(\ti\omega(\Lambda-\lambda,z)-v_\epsilon(\lambda,z)\big)+\tfr12\zeta_1(\lambda)\tr\chi_{\ti\omega}(\Lambda-\lambda,z)\\
	&\hp{=}+\zeta_2(\lambda)\partial_\lambda v_\epsilon(\lambda,z)\\
&=\zeta_1'(\lambda)\br{\int\eta_\epsilon(\lambda-u)\br{\ti\omega(\Lambda-\lambda,z)-\ti\omega(\Lambda-u,z)}du}\\
&\hp{=}+\tfr 12\zeta_2(\lambda)\int\eta_\epsilon(\lambda-u)\big(\tr\chi_{\ti\omega}(\Lambda-u,z)-\tr\chi_{\ti\omega}(\Lambda-\lambda,z)\big)du\\
&\hp{=}+\tfr12\tr\chi_{\ti\omega}(\Lambda-\lambda,z).
}
We conclude 
\eq{\partial_\lambda\ti\omega_\epsilon(\lambda,z)\to \tfr12\tr\chi_{\ti\omega}(\Lambda-\lambda,z)} uniformly as $\epsilon\to0$. Moreover, 
\eq{\tfr12\tr\chi_{\ti\omega}(\Lambda-\lambda,z)>2c} giving $\partial_\lambda\ti\omega_{\epsilon}>c$ for some $\epsilon<\min(\epsilon_0,\frac12\delta)$. A similar argument shows $\partial_\lambda\ti\omega_\epsilon\to 1$ uniformly as $\epsilon\to 0$ on the set $\mathcal{U}_+$. We conclude with some $\epsilon_1$, $\omega_{\epsilon_1}$ satisfies
\eq{\partial_\lambda\ti\omega_{\epsilon_1}(\lambda,z)>0,\q \tr\chi_{\ti\omega_{\epsilon_1}}(\lambda,z)>0}
on $\Omega$, $\ti\omega_{\epsilon_1}(\lambda,z)= \ti\omega(\Lambda-\lambda,z)$ in a neighborhood outside the MOTS of convergence for the mean curvature flow, and $\ti\omega_{\epsilon_1}(\lambda,z)= \lambda$ near $\La_+$.

\begin{defn}
	We say a cross-section $\Si_{out}\sub\Om$ is a strict outermost MOTS if $\Si_{out}:=\{s=\om_{out}(z)\}$ is a MOTS, and given any MOTS, $\Si:=\{s=\om(z)\}$, it follows that $\om_{out}\geq\om$.
\end{defn} 

\subsection{Completion of the proof of \texorpdfstring{\Cref{thm:foliation}}{Theorem 1.4}}
Having established the existence of a global foliation by future un-trapped surfaces, we again for convenience re-parametrize $\Om$ so that the region $$\bigcup_{s\geq 0}\Si_s\sub \Om$$
consists of future un-trapped leaves $\Si_s$, whereby $\Si_{out}:=\Si_0$ corresponds to the MOTS to which the mean curvature flow converges. From \eqref{chi}, we recall the quasi-linear elliptic equation for $\tr\chi_\om$ associated to any graph $\Si_\om:=\{s=\om(z)\}$. From the fact that $\tr\chi_s>0$ for all $s> 0$, we conclude by the strong maximum principle that no graph $\Si_\om$ can satisfy both $\om(p) = \max{\om}>0$ and $\tr\chi_\om(p)=0$. Consequently, $\Si_{out}$ is an outermost MOTS.

\section*{Acknowledgments}
This work was made possible through a research scholarship {JS} received from the DFG and which was carried out at Columbia University in New York. JS would like to thank the DFG, Columbia University and especially Prof.~Simon Brendle for their support. HR would like to acknowledge the support of the National Science Foundation under award No. 1703184.

\bibliographystyle{shamsplain}
\bibliography{Bibliography.bib}
\end{document}